\documentclass[preprint,12pt, twoside]{elsarticle}

\usepackage{amssymb}
\usepackage[utf8]{inputenc}
\usepackage{amsmath}
\usepackage{amsthm}
\usepackage{url}
\usepackage{mathrsfs}
\usepackage{amsfonts}
\usepackage{psfrag}
\usepackage{dsfont}
\usepackage{pifont}
\usepackage{ulem}
\usepackage[dvipsnames,table,xcdraw]{xcolor}
\usepackage{appendix}
\usepackage{mathrsfs}
\usepackage{hyperref}
\usepackage{algorithm, algorithmic}
\usepackage{tikz}
\usepackage{pgfplots}
\pgfplotsset{compat=1.8}
\usepackage{pgfplotstable}
\usepackage{geometry}
\usepackage{todonotes}

\usepgfplotslibrary{groupplots}
\DeclareMathOperator*{\argmin}{\arg\!\min}

\definecolor{PineGreen}{cmyk}{0.92,0,0.59,0.25} 

\journal{Mathematics and Computers in Simulation}

\begin{document}

\begin{frontmatter}

\title{{An SDP dual relaxation for the Robust Shortest Path Problem with ellipsoidal uncertainty: Pierra's decomposition method and a new primal Frank-Wolfe-type heuristics for duality gap evaluation}}

\author[femto,LMB]{Chifaa Al Dahik}
\ead{chifaa.dahik@femto-st.fr}

\author[femto]{Zeina Al Masry}
\ead{zeina.al.masry@ens2m.fr}

\author[LYON]{Stéphane Chrétien}
\ead{stephane.chretien@univ-lyon2.fr}

\author[femto]{Jean-Marc Nicod}
\ead{jean-marc.nicod@ens2m.fr}

\author[LMB]{Landy Rabehasaina}
\ead{landy.rabehasaina@univ-fcomte.fr}

\address[femto]{FEMTO-ST institute, Univ. Bourgogne Franche-Comté, CNRS, ENSMM, Besan\c{c}on, France}
\address[LMB]{Laboratoire de Mathématiques de Besançon, Univ. Bourgogne Franche-Comté, CNRS, Besan\c{c}on, France}
\address[LYON]{Laboratoire ERIC, UFR ASSP, Université Lyon2, Lyon, France}
\date{October 15, 2021}

\begin{abstract}

This work addresses the Robust counterpart of the Shortest Path Problem (RSPP) with a correlated uncertainty set. Since this problem is hard, a heuristic approach, based on Frank-Wolfe's algorithm named Discrete Frank-Wolf (DFW), has recently been proposed. The aim of this paper is to propose a semi-definite programming relaxation for the RSPP that provides a lower bound to validate approaches such as DFW Algorithm. The relaxed problem results from a bidualization that is done {through} a reformulation of the RSPP into a quadratic problem. Then the relaxed problem is solved using a sparse version of Pierra's decomposition in a product space method. This validation method is suitable for large size problems. The numerical experiments show that the gap between the solutions obtained with the relaxed and the heuristic approaches is relatively small.


\end{abstract}

\end{frontmatter}

%

\section{Introduction}

{Robust combinatorial optimization consists in taking uncertainty into account in combinatorial optimization problems. For instance, the robust shortest path problem is the problem of finding the shortest route from a place to another, while the distance (either in terms of time or space) of the different parts of the road are uncertain.} Many definitions of robustness have been proposed in the literature in the context of optimization. The three most common definitions in the context of combinatorial optimization have been formalized in~\cite{kouvelis_robust_2011}. These are absolute robust solution, robust deviation  and relative robust solution. In all these cases, worst case behavior is considered. {For these three definitions, an uncertainty set has to be defined. Many uncertainty sets exist, such as the interval uncertainty, discrete uncertainty, and ellipsoidal uncertainty~\cite{li2011comparative}}. Another family of definitions {is} scenario dependent. In these methods, a decision is taken conditional on the current scenario and the overall optimization problem boils down to a robust two-stage problem~\cite{ben2004adjustable}. This family splits into the notions of  K-adaptability~\cite{hanasusanto2016k}, adjustable robustness, bulk robustness and recoverable robustness. In the case where the data can be considered as governed by a certain probability distribution with unknown parameters, distributionally robust optimization~\cite{rahimian2019distributionally} is also an interesting approach. It consists in choosing the distribution that is most suitable given a robustness criterion. Yet another approach is the notion of almost robust solution~\cite{baron2019almost} that is feasible under most of the realizations and that can use full, partial or no probabilistic information about the uncertain data. Let us also mention that other alternative generic approaches have  also been proposed in the literature: in ~\cite{buhmann2018robust}, a near-optimum solution for several scenarios. Another way to tackle uncertainty that is different from robust optimization is online optimization~\cite{hazan2016introduction}, where decisions are made iteratively, and at each iteration, the problem inputs are unknown, but the decision maker learns from the previous configuration before making his decision. After a decision is made, it is then assessed against the optimal one. Finally, let us add that uncertainty theory was used in another line of work as in~\cite{liu2009some}, for instance. This theory has also been implemented in~\cite{gao2011shortest} in order to give what they call an uncertainty distribution in the case of the shortest path problem.

{This work considers} the absolute robust decision with the uncertainty in the cost function modelled by an ellipsoidal uncertainty set. {The choice of the ellipsoidal uncertainty set is motivated as follows. Unlike the interval uncertainty set, it takes the correlation of the uncertain variables into account, it reduces the combinatorial aspect of the discrete set, it allows the user to control the level of risk that he is ready to take in order to have the right cost. Finally, it leads to a smooth form for the min-max formulation, as shown in Section~\ref{subsec:prob_statement}. This smooth form is well known in portfolio optimization, and it is called mean-risk optimization~\cite{markowitz1952portfolio}}. {In~\cite{kouvelis_robust_2011}}, it is demonstrated that the robust counterparts of easy problems are usually hard to solve, especially if the uncertainty set is not an interval, but is described by an ellipsoidal confidence region. In this case, the robust counterparts of even linear problems become non-linear. To solve these NP-hard problems, {methods exist for the case of non-correlated variables, i.e., for axis-parallel ellipsoids~\cite{poss2018robust, nikolova2010approximation, baumann2014lagrangean, atamturk2008polymatroids}}. {In the case of correlated variables,} branch-and-bound methods exist, as well as improvements by better node relaxations~\cite{buchheim2018robust}. A heuristic approach called DFW (Discrete Frank-Wolfe) for robust optimization {under correlated ellipsoidal uncertainty} based on Frank-Wolfe's algorithm has been proposed in~\cite{al2020frank}. To the best of our knowledge, it is the first algorithm for robust optimization in the ellipsoidal uncertainty adapted for large size problems. 

In order to validate heuristic approaches, one can compare to other exact or heuristic methods, give sub-optimality proofs, or compute lower/upper bounds depending whether it is a minimization or a maximization problem. For minimization problems, lower bounds can be obtained using relaxation schemes such as the ones obtained using Lagrangian dualizations~\cite{boyd2004convex} which often result in solving Semi-Definite Programming (SDP) Problems. 

Semi-Definite Programming (SDP) is a particular class of convex optimization problems which appears in various engineering motivated problems, including the most efficient relaxations of some NP-hard problems such as often encountered in combinatorial optimization or Mixed Integer Programming~\cite{anjos2011handbook}. SDP can be written as minimization over symmetric (resp. Hermitian) positive semi-definite matrix variables, with linear cost function and affine constraints, i.e., problems of the form:
\begin{align}
\min_{Z \succeq 0} \left( \langle A,Z\rangle:  \langle B_j,Z\rangle =b_j \mbox{ for } j=1,\ldots,m\right)
\end{align}
where $A, B_1, \ldots, B_m$ are given matrices. Compact SDPs can be solved in polynomial time. SDP was extensively studied over the last three decades since its early use which  can be traced back to~\cite{scobey1978vector} and~\cite{fletcher1981nonlinear}. In particular, Linear Matrix Inequalities (LMI) and their numerous applications in control theory, system identification and signal processing have been a central drive for the use of SDP in the 90's as reflected in the book~\cite{boyd1994linear}. One of the most influential paper for that era, is~\cite{goemans1995improved} in which SDP was shown to provide a 0.87 approximation to the \textit{Max-Cut} problem, a famous clustering problem on graphs. Other SDP schemes for approximating hard combinatorial problems have subsequently been devised for the graph coloring problem~\cite{karger1998approximate}, for satisfiability problem~\cite{goemans1995improved,goemans1994new}. These results were later surveyed in~\cite{lemarechal1999semidefinite,goemans1997semidefinite} and~\cite{wolkowicz1999semidefinite}. Numerical methods for solving SDP's are manifold and various schemes have been devised for specific structures of the constraints. One of these families of methods is the class of interior point methods~\cite{nesterov1994interior}. Such methods are known to be of the most accurate type, but suffer from being not scalable in practice. Another family of methods is based around the alternating direction method of multipliers (ADMM) technique~\cite{boyd2011distributed}. ADMM approaches are usually faster as they can be implemented in a distributed architecture. As such, they often appear to be faster and more scalable than interior point methods at the price of a worse accuracy. Other methods can also be put to work as the method of Pierra~\cite{pierra1984decomposition} upon which {the present work further elaborates}. 

In this work, the quality of the solution of the DFW heuristic approach is evaluated by computing a lower bound. This lower bound is obtained by a bidualization of the robust problem, that is a SDP relaxation. In order to solve the corresponding SDP problem, the applied method is the decomposition through formalization into a product space proposed in~\cite{pierra1984decomposition}, with a sparse computation to reduce the memory storage necessity. {It is shown} that this algorithm is a validation method that is also adapted for large size problems.

The paper is organized as follows: Section~\ref{sec:robust_SP} presents the robust shortest path problem, and recalls two approaches to solve the problem; in particular we oppose the classical CPLEX exact solving approach to an efficient heuristic algorithm, the so-called Discrete Frank Wolfe (DFW) algorithm, proposed in a previous paper~\cite{al2020frank} that performs well on simulations and is scalable. Since the exact approach is costly, the main contribution of this paper is to propose a validation method for DFW by an efficient relaxation method (SDP) that provides a lower bound for the cost function. This is described in Section~\ref{sec:validation}.
Finally, Section~\ref{sec:experiments} numerically validates this approach by showing that the corresponding gap between the solutions obtained with the relaxed and {the proposed} heuristic approaches is relatively small.

\medskip
\textbf{Notation.} Throughout the paper, the following matrix notations {have been used}. Unless stated otherwise, all vectors belonging to $\mathbb{R}^l$ for some $l\in \mathbb{N}^*$ are column vectors. Furthermore, for some matrix $M$, $M[a:b, c:d]$ denotes, for all integers $a\le b$ and $c\le d $, the sub-block containing the entries in the rows $a$ to $b$ and columns $c$ to $d$. $M[a, c:d]$ (resp. $M[a:b, c]$) is short for $M[a:a, c:d]$ (resp. $M[a:b, c:c]$). $M^T$ is the transpose of $M$. Finally, $0_{(l,l)}$ is the block of dimension $l \times l$ with zeros everywhere and $I_l$ is the identity block of dimension $l \times l$.

%

\section{The robust shortest path problem}
\label{sec:robust_SP}

In this section, the robust shortest path problem with the ellipsoidal uncertainty set is stated. The problem form to solve is then given in order to propose a robust solution. Both an exact and a heuristic method for solving this problem {are presented}.

\subsection{Problem statement}
\label{subsec:prob_statement}
Consider the linear programming form of the shortest path problem~\cite{schrijver2003course}, that can be written as
\begin{align}
\label{P}
\min_{x \in X} c^T x,
\end{align}
with $c \in \mathbb{R}^m$ being the cost vector and $X=\{x \in \{0,1\}^m; Ax=b \}$, where $A\in \mathbb{R}^{n\times m}$ is the so-called incidence matrix corresponding to the underlying graph and $b\in \mathbb{R}^n$ is the vector with $m$ vertices that defines the source and the destination nodes.

{This paper considers} the particular situation where the cost vector $c \in \mathbb{R}^m$  is uncertain, i.e., lies in an uncertainty set $U\subseteq \mathbb{R}^m$. Then the robust counterpart of Problem~\eqref{P} is the following
\begin{align}
\label{RP}
\min_{x \in X} \max_{c \in U} c^T x .
\end{align}
In the particular case where the cost vector $c$ is a random vector following a multinormal distribution with expectation $\mu \in \mathbb{R}^m$ and covariance matrix $\Sigma\in \mathbb{R}^{m\times m}$, then $c$ belongs to the confidence set $E$ with probability $1 - \alpha \in [0,1]$, where $E$ is the following ellipsoid

\begin{align}
\label{E}
E= \{ c \in \mathbb{R}^m ; (c-\mu)^T \Sigma^{-1} (c-\mu) \leq \Omega^2 \},
\end{align}
with $\Omega=\Omega_\alpha\geq 0$ being a function of $\alpha$ that represents the level of confidence.\\

One interesting fact about the ellipsoidal uncertainty set is that the $\min$-$\max$ problem~\eqref{RP} can be reduced to a non-linear programming problem (for more details see Section 2.2.1.1 of~\cite{ilyina2017combinatorial}):
\begin{align}
\label{elli_prob_omega}
\min_{x \in X} \max_{c \in E} c^T x = \min_{x \in X} \mu^T x + \Omega \sqrt{ x^T \Sigma x}.
\end{align}

Without loss of generality, {it is assumed} throughout this paper that $\Omega=1$. This amounts to replace $\Omega^2 \Sigma$ by $\Sigma$, so that Problem~\eqref{elli_prob_omega} reads
\begin{align}
\label{elli_prob}
\min_{x \in X} \mu^T x + \sqrt{ x^T \Sigma x} \ = \ \min_{x \in X} g(x),
\end{align}
where $g(x) = \mu^T x + \sqrt{ x^T \Sigma x}$.\\

 The remaining part of this section addresses the robust shortest path problem by solving Problem~\eqref{elli_prob}. This problem is a non-linear non-convex problem{, so it is} challenging to find an appropriate method to solve it.

\subsection{Exact method for solving the robust problem}

In order to solve the robust shortest path problem, {solving} Problem~\eqref{elli_prob} {is needed}. One possible way is to solve this problem in two steps. First, rewrite it as a Binary Second Order Cone Programming problem (BSOCP) (Problem~\eqref{cone}). The solution is obtained by solving Problem~\eqref{cone} in the second step. Problem~\eqref{cone} {is stated} as follows:

\begin{align}
\label{cone} 
\min &\quad\mu^T x+z \\
\text{s.t. } &\quad (y,z)^T \in K_{m+1} \notag \\
&\quad y=(\Sigma^{\frac{1}{2}})^T x \notag \\
&\quad x \in X , y \in \mathbb{R}^m, z \in \mathbb{R}_+ \notag,
\end{align}
with
$K_{m+1}=\{x \in \mathbb{R}^{m+1}; \| (x_1, \ldots,x_m)^T\|_2 \leq
x_{m+1}\}$ being a \textit{second order cone}. The calculations are detailed in~\cite[Section 4]{al2020frank}.

Problem~\eqref{cone} can be solved by branch-and-bound methods, and existing BSOCP solvers such as CPLEX~\cite{ibm_2019} solve this problem. However, for large size problems, it is no longer possible to use branch-and-bound methods, since their time complexity is exponential, and in their worst case, they may need to explore all the possible permutations of the combinatorial problem at hand.
Thus a heuristic algorithm named DFW Algorithm (Discrete Frank-Wolfe) has been proposed in~\cite{al2020frank} and is presented in the following section.

\subsection{A heuristic approach based on Frank-Wolfe}

As proved in the previous section, solving the robust shortest path problem with a scalable heuristic approach seems mandatory for large size problems. The heuristic algorithm proposed in~\cite{al2020frank} to solve Problem~\eqref{elli_prob} is based on the Frank-Wolfe algorithm~\cite{frank1956algorithm}. On one side, the classical Frank-Wolfe (FW) algorithm is a convex optimization algorithm that proceeds by moving towards a minimizer of the linear appromixation of the function to minimize. The heuristic DFW algorithm in turn uses the classical Frank-Wolfe algorithm to minimize $g(x) =\mu^T x + \sqrt{ x^T \Sigma x}$ over the convex hull of $X$, and {due} to the {integrality of the relaxation}, the intermediate gradient steps are feasible solutions for the discrete problem. These gradient steps are good feasible solutions. DFW Algorithm returns the best of these intermediate steps as an approximate solution: more concretely, it is the one that minimizes the objective function $g$ among the discovered feasible solutions. DFW Algorithm is detailed in Algorithm~\ref{algo}.

\begin{algorithm}[hbt!]
\caption{DFW: a Frank-Wolfe based algorithm to solve (\ref{elli_prob})}
\begin{algorithmic}[1]

  \STATE $x^{(0)}$ a random feasible solution, $\varepsilon>0$ close to zero, $K$ the maximum number of iterations.
  \STATE $k \leftarrow 1$
  \STATE stop $\leftarrow$ false
  \WHILE {$k\leq K$ and $\neg$stop}
  \IF {$g(x^{(k-1)})-g(x^{(k)})<\varepsilon$} \STATE stop $\leftarrow$ true
  \ELSE
  \STATE Choose $\displaystyle s^{(k)} \in \argmin_{y \in \text{Conv}(X)} \nabla g(x^{(k)})^T y \text{, with } s^{(k)} \in X$
  \STATE $\displaystyle \gamma^{(k)} \leftarrow \argmin_{\alpha \in [0,1]} g(x^{(k)}+\alpha(s^{(k)}-x^{(k)}))$
  \STATE $x^{(k+1)} \leftarrow x^{(k)}+\gamma^{(k)}(s^{(k)}-x^{(k)})$
  \ENDIF
  \STATE $k++$
  
  \ENDWHILE

  \STATE \textbf{return} $\displaystyle \argmin_{s \in \{s^{(1)}, \ldots, s^{(k-1)}\}} g(s)$

\end{algorithmic}
\label{algo}
\end{algorithm}

Note $x^*$ the optimal solution of Problem~\eqref{elli_prob}, and $\hat{x}$ the approximate solution given by DFW Algorithm. The aim of the next section is to evaluate the quality of the solution $\hat{x}$.

%

\section{A lower bound by SDP relaxation}
\label{sec:validation}

The first way to evaluate the quality of the solution given by the DFW Algorithm or any other approach that solves Problem~\eqref{elli_prob} is to compare it with the solution given by the optimal solution of the BSOCP using an exact solver like CPLEX (see the previous section). Since this approach is no longer usable when considering large size problem, an other option to evaluate the quality of the solution has to be proposed. To do so a lower bound by bidualization with an efficient algorithm to solve the corresponding problem is presented in this section.

\subsection{Bidualization of a quadratic problem} 
\label{bidualsection}

Before giving a lower bound for Problem~\eqref{elli_prob}, a lower bound by bidualization for any quadratic problem {is stated}. Then, Problem~\eqref{elli_prob} {is written} as a quadratic problem following the general form.

A lower bound for Quadratic Programming problems is proposed in~\cite{lemarechal1999semidefinite}. This lower bound is the solution of a bidual problem that is written in the form of a Semi-Definite Programming (SDP) problem. This bidualization procedure is nothing but the very known SDP relaxation, as in the following.

Consider the quadratic problem with $N$ constraints:
\begin{align}
\label{QP}
\inf q_0(x), \quad x \in \mathbb{R}^d \\
q_j(x)=0, \quad j=1, \ldots,N ,\notag
\end{align}
where
\begin{equation}
q_j(x)=x^T Q_j x + b_j^T x + c_j,\quad j=0, \ldots N \notag
\end{equation}
are $N+1$ quadratic functions defined on $\mathbb{R}^d$, $d\in \mathbb{N}^*$ being the dimension of the problem, with the matrices $Q_j$ lying in the set $S_d$ of symmetric matrices of size $d\times d$, the values $b_j$ in $\mathbb{R}^d$, and the values $c_j$ in $\mathbb{R}$ $ \forall j \in\{ 1, 2, \ldots, N\}$; it is assumed that $c_0=0$.

Applying Lagrangian duality on Problem~\eqref{QP}, and applying duality again reproduces the bidual problem of~\eqref{QP} that is given by
\begin{align}
\label{bidual_QP}
\inf Q_0 \bullet X + b_0^T x, \quad X \in S_d, x \in \mathbb{R}^d, \\ 
Q_j \bullet X + b_j^T x + c_j=0, \quad j=1, \ldots, N, \notag \\
\left[
\begin{array}{c c}
1 & x^T \\
x & X
\end{array}
\right] \succeq 0 ,\notag
\end{align}
where the inner product between the matrices $A$ and $B$ of size $d\times d$ is defined by
\begin{equation}\label{def_innerprod}
A \bullet B =tr(A^T B),
\end{equation}
and where the notation $M\succeq 0$ means that $M$ is positive semi-definite, for any symmetric matrix $M$.

This bidualization has another interpretation: it is also a direct convexification of Problem~\eqref{QP}. Indeed, noticing that~\eqref{QP} can also be written as
\begin{align}
\label{QP_another_writing}
\inf Q_0 \bullet X + b_0^T x, \quad X \in S_d, x \in \mathbb{R}^d, \\ 
Q_j \bullet X + b_j^T x + c_j=0, \quad j=1, \ldots, N, \notag \\
X=x x^T \notag
\end{align}
by setting $X=x x^T$, and writing a quadratic form $x^T Q x$ as $Q \bullet  x x^T$, and then relaxing the nonconvex constraint $X=x x^T$ to $X \succeq x x^T$, that is convex with respect to $(x,X)$. Then, the previous bidualization {can be seen} as a convexification.

Thus, if $p^*$ is the optimal solution of~\eqref{QP}, and $d^{**}$ is the optimal solution of~\eqref{bidual_QP}, then the following inequality holds (see~\cite[Proposition 4.5]{lemarechal1999semidefinite}):

\begin{equation}
\label{inequality}
d^{**} \leq p^*.
\end{equation}

Hence, solving the SDP problem~\eqref{bidual_QP} enables to obtain a lower bound for $p^*$. In general, {this technique is used for the validation of a heuristic method without comparison with the optimal solution}. In this case, {Problem}~\eqref{bidual_QP} is easier, since it is a convex problem. Solving~\eqref{bidual_QP} gives the lower bound $d^{**}$. The distance between the lower bound and the heuristic solution indicates how far this heuristic solution is from the optimal solution. {In other research directions, lower bound can be coupled with a branch-and-bound algorithm for computing an optimal solution. However, the focus in the present paper is on  proposing  a  much  cheaper heuristics  than  the  branch-and-bound  approach,  namely  the  DFW  method. In  order  to validate  this  heuristics, the  quality  of  the  obtained  primal  solution  using  a lower bound obtained by solving a polynomial time SDP problem has been evaluated.}

\subsection{Using the bidualization to compute a lower bound}

\subsubsection{Bidualization of the addressed problem}

{This section aims to} show how to use the bidualization, explained in Section~\ref{bidualsection}, to compute a lower bound for~\eqref{elli_prob}. Recall that Problem~\eqref{elli_prob} has another formulation, that is a BSOCP (Problem~\eqref{cone}). Rewriting~\eqref{cone} more explicitly gives
\begin{align}
\label{cone_explicit} 
\min &\quad \mu^T x+z \\
\text{s.t. } &\quad \sqrt{y^T y} \leq z \notag \\
&\quad y=(\Sigma^{\frac{1}{2}})^T x \notag \\
&\quad x \in X , y \in \mathbb{R}^m, z \in \mathbb{R}_+ .\notag
\end{align}

First, the BSOCP formulation~\eqref{cone_explicit} of~\eqref{elli_prob} can be written as a Binary Quadratic Problem (BQP) since the variables $y$ and $z$ in~\eqref{cone_explicit} are such that $\sqrt{y^T y} \geq 0$ and $z \geq 0$ for any $y \in \mathbb{R}^m$ and $z \in \mathbb{R}_+ $. Thus, Problem~\eqref{cone_explicit} is equivalent to
\begin{align}\label{BQP}
\min &\quad \mu^T x+z \\
\text{s.t. } &\quad y^T y \leq z^2 \notag \\
&\quad y=(\Sigma^{\frac{1}{2}})^T x \notag \\
&\quad x \in X , y \in \mathbb{R}^m, z \in \mathbb{R}_+ .\notag
\end{align}

In order to formulate~\eqref{BQP} as a problem in the form~\eqref{QP}, all the constraints {have to been written} in the form of equalities. {First, the following equivalence holds}
\[
x \in X \iff Ax=b \text{ and } x \in \{0,1\}^m \iff Ax=b \text{ and } x_i(x_i-1)=0 \quad i=1, \ldots, m.
\]

Second, the inequality $y^T y \leq z$ can be transformed into an equality by considering additional variables $c_1$ and $c_2$ as follows:

\[
y^T y \leq z^2 \iff \exists c_1 \in \mathbb{R} \ ; \quad y^T y - z^2 = -c_1^2 \iff  \exists c_1 \in \mathbb{R} \ ; \quad y^T y - z^2 + c_1^2 = 0 $$
$$z \geq 0 \iff \exists c_2 \in \mathbb{R} \ ; \quad z = c_2^2 \iff \exists c_2 \in \mathbb{R} \  ; \quad z - c_2^2 = 0.
\]

The problem~\eqref{BQP} is then equivalent to the following problem
\begin{align}
\label{BQP_explicit} \min &\quad \mu^T x+z \\
\text{s.t. } &\quad  y^T y - z^2 + c_1^2=0 \notag \\
&\quad y=(\Sigma^{\frac{1}{2}})^T x \notag \\
&\quad Ax=b \notag \\
&\quad x_i(x_i-1)=0 \quad i=1, \ldots, m \notag \\
&\quad z - c_2^2= 0 \notag \\
&\quad x \in \mathbb{R}^m, y \in \mathbb{R}^m, z \in \mathbb{R}, c_1 \in \mathbb{R}, c_2 \in \mathbb{R}  \notag.
\end{align}

{Now, Problem~\eqref{BQP_explicit} is written} in a more compact way , i.e., in function of one vector variable $u=[x,y,z,c_1,c_2] \in \mathbb{R}^{2m+3}$, and write each constraint individually. This makes Problem~\eqref{BQP_explicit} equivalent to

\begin{align} \label{BQP_u} 
\min &\quad (\tilde{\mu} + \delta_{2m+1})^T u \\
\text{s.t. } &\quad u^T(\mathds{1}_y^T \mathds{1}_y  - \delta_{2m+1 \text{ } 2m+1} + \delta_{2m+2 \text{ } 2m+2}) u=0 \notag\\
&\quad (\mathds{1}_y - (\tilde{\Sigma^{\frac{1}{2} T}}))_i^T u= 0 \quad i=1, \ldots, m \notag \\
&\quad \tilde{A}_j^T u=b_j \quad j=1, \ldots, n \notag \\
&\quad u^T \delta_{ii} u - \delta_i^T u =0 \quad i=1, \ldots, m \notag \\
&\quad - u^T \delta_{2m+3 \text{ }2m+3} u + \delta_{2m+1}^T u = 0 \notag \\
&\quad u \in \mathbb{R}^{2m+3}  \notag,
\end{align}
where the vectors and matrices that appear in Problem~\eqref{BQP_u} are defined as follows 
\begin{itemize}
\item the vector $\tilde{\mu}$ of size $2m+3$ is defined block-wise as $\tilde{\mu}=[\mu,0, \ldots,0]^T$, so that $\tilde{\mu}^T u= \mu^T x $ if $u=[x,y,z,c_1,c_2]$, 
\item for any $k=1, \ldots, 2m+3 $ , $\delta_k \in \mathbb{R}^{2m+3}$ is such that $\delta_k (l)=1$ if $k=l$, and $0$ if else, so that $\delta_{2m+1}^T u=u_{2m+1}=z$, and $\delta_i^T u=x_i$ for $i=1, \ldots m$,
\item  $\mathds{1}_y$ is an $m \times (2m+3)$ matrix such that $\mathds{1}_y [m+1:2m;m+1:2m]=I_m$ and $0$ elsewhere, so that $\mathds{1}_y u =y$ and $u^T\mathds{1}_y^T \mathds{1}_y u=y^Ty$,
\item for any $i,j=1, \ldots,2m+3$, $\delta_{i,j}$ is a $(2m+3) \times (2m+3)$ matrix, such that $\delta_{i,j} (k,l)=1$ if $i=k$ and $j=l$, and $0$ if else. So that $u^T \delta_{2m+1, 2m+1} u=z^2$, $u^T \delta_{2m+2,2m+2} u=c_1^2$, $u^T \delta_{2m+3, 2m+3} u=c_2^2$ and $u^T \delta_{ii} u=u_i^2$ for $i=1, \ldots, m$,
\item $\tilde{\Sigma^{\frac{1}{2} T}}$ is an $m \times (2m+3)$ matrix such that $\tilde{\Sigma^{\frac{1}{2} T}}[1 : m;1 : m]=\Sigma^{\frac{1}{2} T}$ and the other entries are zeros, so that $\tilde{\Sigma^{\frac{1}{2} T}} u = \Sigma^{\frac{1}{2} T} x $,
\item $\tilde{A}$ is an $n \times (2m+3)$ matrix such that $\tilde{A}[1 : n;1 : m]=A$ and the other entries are zeros, so that $\tilde{A} u = A x $.
\end{itemize}
Then, the bidual problem of~\eqref{BQP_u} is the following
\begin{align}
\label{bidual_u} 
\min &\quad (\tilde{\mu} + \delta_{2m+1})^T u \\
\text{s.t. } &\quad (\mathds{1}_y^T \mathds{1}_y  - \delta_{2m+1 \text{ } 2m+1} + \delta_{2m+2 \text{ } 2m+2}) \bullet U=0 \notag\\
&\quad (\mathds{1}_y - (\tilde{\Sigma^{\frac{1}{2} T}}))_i^T u= 0 \quad i=1, \ldots, m\notag \\
&\quad \tilde{A}_j^T u=b_j \quad j=1, \ldots, n \notag \\
&\quad  \delta_{ii} \bullet U - \delta_i^T u =0 \quad i=1, \ldots, m \notag \\
&\quad -  \delta_{2m+3 \text{ }2m+3} \bullet U + \delta_{2m+1}^T u = 0 \notag \\
&\quad \left[
\begin{array}{c c}
1 & u^T \\
u & U
\end{array}
\right] \succeq 0,  \quad U \in S^{2m+3}, \quad u \in \mathbb{R}^{2m+3}.  \notag
\end{align}
The last step consists in writing~\eqref{bidual_u} in a compact way {with} the {following} change of variable 
\begin{align}
Z = \left[
\begin{array}{c c}
1 & u^T \\
u & U
\end{array}
\right] \in S^{2m+4}. \notag
\end{align}
This can be done using the following changes: 
\begin{enumerate}
\item For any $v \in \mathbb{R}^{2m+3}$, write $v^T u= V \bullet U$,
where $V \in S^{2m+4}$ is defined by
\begin{align}
V={\frac{1}{2}}\left[
\begin{array}{c|c}
0 & v^T\\
\hline
v & 0
\end{array}
\right] \in S^{2m+4}\notag
\end{align}
\item For any $W \in S^{2m+3}$, write $ W \bullet U = \textbf{W} \bullet Z$,
where $\textbf{W}\in S^{2m+4}$ is defined by
\begin{align}
\textbf{W}=\left[
\begin{array}{c|c}
0 & \ldots 0 \\
\hline
0 & W
\end{array}
\right]. \notag
\end{align}

\end{enumerate}
{As a result of} this change of variable, the bidual problem of~\eqref{BQP} can be written as an SDP problem in the more compact way:
\begin{align}
\label{bidual_y}
\min &\quad M \bullet Z \\ 
\text{s.t.} &\quad Z \in S^{2m+4}\notag \\
&\quad O_j \bullet Z = b_j, j=1, \ldots, n ,\notag \\
&\quad C_i \bullet Z = 0, i=1, \ldots, m ,\notag \\
&\quad Q \bullet Z = 0, \notag \\ 
&\quad D_i\bullet Z =0, \quad i=2, \ldots, m+1 \notag \\
&\quad L \bullet Z = 0, \notag \\
&\quad Z \succeq 0, \notag 
\end{align}
where $M\in S^{2m+4}$ is defined as follows
\begin{align}
M={\frac{1}{2}}\left[
\begin{array}{c|c}
0 & \mu^T 0\ldots 0\,1\,0\,0\\
\hline
\mu & \\
0 & \\
\vdots \\
0 & \mbox{\huge $0$} \\
1 &\\
0 &\\
0 &
\end{array}
\right], \notag
\end{align}
that is $M[1,2:m+1]=\frac{1}{2} \mu^T$, $M[1,2m+2]=\frac{1}{2}$,
$M[2:m+1,1]=\frac{1}{2} \mu$, $M[2m+2,1]=\frac{1}{2}$, and zero elsewhere. The matrix $O_j\in S^{2m+4}$ is defined for all $j=1, \ldots, n$ by
\begin{align}
O_j={\frac{1}{2}}\left[
\begin{array}{c|c}
0 & A_j^T 0\ldots 0\,0\,0\,0\\
\hline
A_j & \\
0 & \\
\vdots \\
0 & \mbox{\huge $0$} \\
0 &\\
0 &\\
0 &
\end{array}
\right], \notag
\end{align}
that is $O_j[1,2:m+1]=\frac{1}{2} A_j^T$,
$O_j[2:m+1,1]=\frac{1}{2} A_j$, and zero elsewhere.
The matrix $C_i\in S^{2m+4}$ is defined for all $i=1, \ldots,m$ by
\begin{align}
C_i={\frac{1}{2}}\left[
\begin{array}{c|c}
0 & -(\Sigma^{\frac{1}{2} T})_i^T 0\ldots\,0\,1\,0\ldots0\\
\hline
-(\Sigma^{\frac{1}{2} T})_i  & \\
0 & \\
\vdots \\
0 & \\
1 & \mbox{\huge $0$} \\
0 &\\
\vdots & \\
0 &
\end{array}
\right], \notag
\end{align}
that is $C_i[1,2:m+1]=-\frac{1}{2}(\Sigma^{\frac{1}{2} T})_i^T$, $C_i[1,m+1+i]=\frac{1}{2}$, $C_i[2:m+1,1]=-\frac{1}{2}(\Sigma^{\frac{1}{2} T})_i$, $C_i[m+1+i,1]=\frac{1}{2}$, and zero elsewhere. {Also,} define the matrix $Q$ by
\begin{align}
Q=\left[
\begin{array}{c|c | c | c c c}
0 & 0 \ldots 0 & 0 \ldots 0 & 0 & 0 & 0\\
\hline
 & 0_{(m,m)} & & & &\\
\hline
 & & I_m & & &\\
\hline
\vdots & & & -1 & & \\
 & & & & 1 & \\
0 & & & & & 0
\end{array}
\right], \notag
\end{align}
that is $Q[m+2:2m+1,m+2:2m+1]$ is the identity matrix of dimension $m$, $Q[2m+2,2m+2]=-1$, $Q[2m+3,2m+3]=1$  and zero elsewhere.
Next, for the definition of the matrices $D_i$, for every $i=2,m+1$, $D_i$ is a $2m+4 \times 2m+4$ matrix such that $D_i[i,i]=1$, $D_i[i,1]=-\frac{1}{2}$ and $D_i[1,i]=-\frac{1}{2}$.
Finally, $L$ is a $2m+4 \times 2m+4$ matrix such that $L[1, 2m+2]=\frac{1}{2}$, $L[2m+2,1]=\frac{1}{2}$ and $L[2m+4,2m+4]=-1$.

\subsubsection{The biduality gap}
Now that the bidual problem~\eqref{bidual_y} of~\eqref{BQP} {is stated}, the lower bound inequality~\eqref{inequality} reads here
\begin{equation}
\text{val}(\eqref{bidual_y})\leq \text{val}(\eqref{BQP}) ,\notag
\end{equation}
where $\text{val}((P))$ denotes the optimal value for a given problem (P). {As a result of} the equivalence between~\eqref{elli_prob} and~\eqref{BQP}, val(\eqref{BQP}) equals $g(x^*)$. This gives us an additional inequality:
\begin{equation}
\text{val}(\eqref{bidual_y})\leq \text{val}(\eqref{BQP}) \ = g(x^*) \leq g(\hat{x}) .\notag
\end{equation}
Or, written differently,
\begin{equation}
\label{the_used_inequality}
d^{**} \leq p^*= g(x^*) \leq g(\hat{x}).
\end{equation}
Thus, $d^{**}$ is a lower bound that allows to evaluate the quality of the heuristic solution of the DFW Algorithm. Hence, the biduality gap $BG$ is defined as
\begin{equation}
BG= g(\hat{x}) -d^{**}.
\label{bg}
\end{equation}
A corresponding relative gap $RBG$ is defined as
\begin{equation}
RBG=\frac{g(\hat{x})-d^{**}}{d^{**}}.
\label{rbg}
\end{equation}
More explicitly, the validation process is the following: first solve {the} robust shortest path problem using the heuristic approach DFW and find a heuristic solution $\hat{x}$. Then evaluate the quality of this solution using Inequality~\eqref{the_used_inequality} by proceeding as follows. {Then, $d^{**}$ is computed}, and if the gap between $d^{**}$ and $g(\hat{x})$ is small, then the gap between $g(x^*)$ and $g(\hat{x})$ is small too, since $g(\hat{x}) - d^{**} \geq g(\hat{x}) - g(x^*) \geq 0$. The only missing step now is to compute $d^{**}$. {The next section shows} how to solve~\eqref{bidual_y} to compute $d^{**}$.

\subsection{Solving the SDP problem} \label{subsection:SDPsolving}

The above sections aim at showing that a lower bound for the robust shortest path problem is the solution of an SDP problem that has to be solved. As detailed in the introduction, interior point methods are used to solve SDP problems, which gives a first way to solve the SDP problem~\eqref{bidual_y}: an option is to implement this resolution using the CVXPY Python package~\cite{diamond2016cvxpy} which is a Python-embedded modeling language for convex optimization problems.
CVXPY converts the convex problems into a standard form known as conic form, a generalization of a linear program.  The conversion is done using graph implementations of convex functions. The resulting cone program is equivalent to the original problem, so {solving it gives} a solution of the original problem. In particular, it solves the semi-definite programs using interior point methods. It is rather simple to use CVXPY to solve {the} SDP problem~\eqref{bidual_y}: define the function to minimize, the constraints of the problem, and then launch the solver. However this simplicity has a price: the problem definition requires the storage of the matrices that describe the problem. More precisely, {there are} $n+2m+4$ matrices of dimension $2m+4 \times 2m+4$: one matrix to define the objective function, and $n+2m+3$ matrices for the constraints. This is a significant issue because of the storage necessity, especially in large size problems. To illustrate how big the storage grows with respect to the problem size, take a medium grid graph with $10 \times 10$ nodes ($n=100$, $m=360$). This problem size requires the storage of $824$ matrices of dimension $724 \times 724$ (this takes 3.45 Gigabytes in double precision). A relatively big grid graph with $40 \times 40$ nodes ($n=1600$, $m=6240$) needs the storage of $14084$ matrices of dimension $12484 \times 12484$ (this takes 17.5 terabytes in double precision). But since most of the matrices are sparse, another efficient approach {is proposed,} where sparse computations {are done}, {which} allows us to avoid this main drawback considering the matrices storage. Before tackling this memory storage issue, {the following} describes the practical algorithm that {has been} implemented in order to find this $d^{**}$.

\subsubsection{Pierra's Decomposition through formalization in a product space}

Consider a general minimization problem in a finite dimensional Hilbert space $\mathcal{H}$ equipped with a norm $\|.\|_2$. Suppose {that the goal is to solve the problem}
\begin{align}
\label{min_general}
\min_{x \in \mathcal H} & \quad f(x) \\
\text{ s.t. } & \quad x \in \cap_{j=1}^J \mathcal{S}_j ,\notag
\end{align}
where $f$ is a differentiable function, and $\mathcal{S}_1, \ldots, \mathcal{S}_J$ are convex subsets of ${\mathcal H}$. Exploiting the fact that the constraint space is an intersection of convex sets, Pierra in~\cite{pierra1984decomposition} proposes a method for solving Problem~\eqref{min_general}. This is described in Algorithm~\ref{Pierra_general}, where, for a function $h:\mathcal{H}\longrightarrow \mathbb{R}$, the proximal function associated to $h$, which is defined in~\cite[Theorem 3.2]{pierra1984decomposition}, is given by

\begin{equation}\label{def_prox_h}
Prox_h(y)=\argmin_{x\in \mathcal{H}}\left[h(x)+\frac{1}{2} \|x-y\|_2^2 \right],
\end{equation}
and where $I_{\mathcal{S}_j}(x)$, referring to the indicator function for the set $\mathcal{S}_j$, equals $0$ if $x\in \mathcal{S}$ and $+\infty$ otherwise. Finally, $\varepsilon>0$ is a tuning parameter for the minimization step which value is small (e.g. $\varepsilon = 10^{-4}$).

\begin{algorithm}[hbt!]
\caption{Pierra's algorithm to solve~\eqref{min_general}}
\begin{algorithmic}[1]

    \STATE $x^0 \in \mathcal{H}$ random, $k \in \mathbb{N}$, $\lambda \in ]0,1]$, $\varepsilon$ small, P the maximum number of iterations.
    \STATE $p \leftarrow 0$
    \STATE stop $\leftarrow$ false
    \WHILE {$p\leq P$ and $\neg$stop}
    \STATE  $v_j^{p+1} \leftarrow Prox_{I_{\mathcal{S}_j}+\frac{1}{2J}\varepsilon f} (x^p), j=1, \ldots, J$
        \STATE $b^{'p+1} \leftarrow \frac{1}{J} \Sigma_{j=1}^{J} v_j^{p+1}$
	    \IF {$b'^{p+1}=x^p$} 
	        \STATE stop $\leftarrow$ true
        \ELSE
	        \STATE $b^{p+1} \leftarrow x^p+\beta^{p+1} (b'^{p+1}-x^p)$ with $\beta^{p+1} \leftarrow \frac{\Sigma_{j=1}^{J} \|v_j^{p+1}-x^p\|^2}{J\|b'^{p+1}-x^p\|^2}$
            \STATE $x^{p+1} \leftarrow 
            \begin{cases}
            x^p+\lambda(b^{p+1}-x^p), & \text{if } p+1 \equiv k (\textrm{mod } k).\\
            b^{p+1} & \text{otherwise}.
            \end{cases}$
		 \STATE $p++$
		 \ENDIF
    \ENDWHILE
    \STATE \textbf{return} $x^p$

\end{algorithmic}
\label{Pierra_general}
\end{algorithm}

The idea of this algorithm comes from the formalization of the constraint set $\cap_{j=1}^J \mathcal{S}_j$ introducing the set $\textbf{H}=\mathcal{H}^J$. Indeed, {defining} $\textbf{S}=\mathcal{S}_1 \times \cdots \times \mathcal{S}_J$, and {denoting} the diagonal convex $\textbf{D}$ as the subspace of $\textbf{H}$ of all the vectors of the form $(x, \ldots,x)$, with $x\in \mathcal{H}$, {implies that} Problem~\eqref{min_general} can be reformulated in $\textbf{H}$ as a minimization problem over $\textbf{S} \cap \textbf{D}$.
To solve~\eqref{min_general}, Pierra's Algorithm can be described in three steps: (i) The first step (line 5 of Algorithm~\ref{Pierra_general}) comes from the projection on $\textbf{S}$, with a part of minimization of the objective function. Here, the proximal function can be explained intuitively as follows: for every constraint space $S_j$, it both minimizes the function $f$ and stays close to $x^p$, and since $x^p$ partially results from a point that belongs to all the constraint spaces, then $x^p$ converges to the optimal solution; (ii) The second step comes from the projection over the diagonal convex $\textbf{D}$, represented in line 6 of Algorithm~\ref{Pierra_general}; (iii) Finally, the third step is an extrapolation step. In simple words, the extrapolation represented in line 11 is used to center the iterate $x^p$ from time to time, at each $k$ iterations: in~\cite[Section 4]{pierra1984decomposition}, it is explained that without the centring technique, the convergence seems to become ineffective, and on the other side, centring at each iteration can lead to an ineffective extrapolation. It has been proved in~\cite[Theorem 3.3]{pierra1984decomposition} that this algorithm converges. All the theoretical background of Pierra's algorithm can be found in~\cite{pierra1984decomposition}.

\subsubsection{Pierra's algorithm adapted to solve {the considered} SDP problem} 

{This part aims} to apply Algorithm~\ref{Pierra_general} to solve~\eqref{bidual_y}.
In this case, the corresponding Hilbert space is set as ${\cal H}=S^{2m+4}$, with the norm $\|.\|_F$ that is associated to the inner product $\bullet$ defined in Section~\ref{bidualsection} by~\eqref{def_innerprod}, such that $\|A\|_F^2=tr(A^T A)$. The function to minimize in Problem~\eqref{min_general} is given by $ f:Z\in S^{2m+4}\mapsto f(Z)=M \bullet Z$, and the integer $J$ equals $n+2m+3$. The convex sets $\mathcal{S}_1, \ldots,\mathcal{S}_J$ are defined as follows:
\begin{equation}\label{form_Sj}
\begin{array}{rcl}
\mathcal{S}_j&= &\{Z \in S^{2m+4};\ \textbf{A}_j \bullet Z =\textbf{b}_j\},\quad j=1, {\ldots}, n+2m+2,\\
 \mathcal{S}_J&=& \{Z \in S^{2m+4};\ Z \succeq 0 \},
\end{array}
\end{equation}
where $\textbf{A}_j$, $\textbf{b}_j$, $j=1, {\ldots}, n+2m+2$ are respectively matrices and scalars defined by
\begin{equation}\label{form_Aj_bj}
\textbf{A}_j=
\left\{
\begin{array}{cl}
 O_j,& j=1, \ldots, n, \\
 C_{j-n},& j=n+1, \ldots, n+m,\\
 Q,& j= n+m+1,\\
 D_{j-(n+m+1)},& j=n+m+2, \ldots, n+2m+1,\\
 L, & j=n+2m+2,
\end{array}
\right.
\quad
\textbf{b}_j = \left\{
\begin{array}{cl}
 b_j,& j=1, \ldots, n, \\
 0, & j=n, \ldots, n+2m+2.
\end{array}
\right.
\end{equation}

In {the considered} case, the proximal function associated to $I_{S_j}+\frac{1}{2J}\varepsilon f$ on line 5 of Algorithm~\ref{Pierra_general} is computed using the definition~\eqref{def_prox_h} in the following way:
\begin{align}
Prox_{I_{\mathcal{S}_j}+\frac{1}{2J}\varepsilon f} (x^p) &=
  \argmin_{Z \in \mathcal{S}_j} \left[\frac{1}{2J} \varepsilon M \bullet Z + \frac{1}{2} \| Z - x^p \|_F^2 \right]   \notag \\
&= \argmin_{Z \in \mathcal{S}_j} \left[\frac{1}{2J} \varepsilon M \bullet Z + \frac{1}{2} \| Z \|_F^2 -Z \bullet x^p + \frac{1}{2} \| x^p \|_F^2 \right]    \notag \\
&= \argmin_{Z \in \mathcal{S}_j} \left[ \frac{1}{2} \| Z \|_F^2  -Z \bullet (x^p -  \frac{1}{2J} \varepsilon M )  + \frac{1}{2} \| x^p \|_F^2 \right]     \notag \\
&= \argmin_{Z \in \mathcal{S}_j}  \left[\frac{1}{2} \| Z  - (x^p -  \frac{1}{2J} \varepsilon M )\|_F^2 = Proj_{\mathcal{S}_j}(x^p-\frac{1}{2J} \varepsilon M)\right], \label{prox}
\end{align}
where $Proj_{\mathcal{S}_j}$ is the projection on the set $\mathcal{S}_j$. Thus, one sees from~\eqref{prox} that there remains to compute the projections over the constraint spaces defined by~\eqref{form_Sj}. Those spaces are of two kinds. First, for any constraint in the form $C=\{Z \in S^{2m+4}; \textbf{A} \bullet Z =\textbf{b}\}$, the following explicit projection formula {holds}:
\begin{equation}
Proj_C(Z)=Z+\left(\frac{\textbf{b}-\textbf{A} \bullet Z}{\|\textbf{A}\|_F^2}\right) \textbf{A} \notag.
\end{equation}
Second, concerning the projection on the constraint space $\mathcal{S}_J=\{Z \in S^{2m+4};\ Z \succeq 0\}$,

\begin{equation}
Proj_{\mathcal{S}_J} (Z)= U \max \{\Lambda,0\} U^T, \notag
\end{equation}
where $Z=U \Lambda U^T$  is the eigenvector decomposition of the matrix $Z$ (see~\cite[section 20.1.1]{anjos2011handbook}).
In view of all these considerations, Pierra's algorithm  applied on problem~\eqref{bidual_y} is described in Algorithm~\ref{algo2}.
\begin{algorithm}[hbt!]
\caption{Pierra's algorithm to solve the SDP problem~\eqref{bidual_y}}
\begin{algorithmic}[1]

  \STATE $Z^1 \in S^{2m+4}$ random, $k \in \mathbb{N}$, $\lambda \in ]0,1]$, $\varepsilon$ small, $\alpha$ small, $P$ the maximum number of iterations.
  \STATE $p \leftarrow 1$
  \STATE stop $\leftarrow$ false
  \WHILE {$p \leq P$ and $\neg$stop}
        \STATE $Y^p \leftarrow Z^p-\frac{\varepsilon}{2(n+2m+3)}M $
		\FOR {$j=1$ to $n$}
		    \STATE $Z_{j}^{p+1} \leftarrow Y^p + (\frac{b_j-O_j \bullet Y^p}{\|O_j\|^2}) O_j$
		\ENDFOR
		\FOR {$i=1$ to $m$}
		    \STATE $Z_{n+i}^{p+1} \leftarrow Y^p + (\frac{-C_i \bullet Y^p}{\|C_i\|^2}) C_i$
		\ENDFOR
		
		\STATE $Z_{n+m+1}^{p+1} \leftarrow Y^p+(\frac{0-Q \bullet Y^p}{\|Q\|^2}) Q$
		\FOR {$i=2$ to $m+1$}
		    \STATE $Z_{n+m+i}^{p+1} \leftarrow Y^p + (\frac{0-D_i \bullet Y^p}{\|D_i\|^2}) D_i$
		\ENDFOR
		\STATE $Z_{n+2m+2}^{p+1} \leftarrow Y^p+(\frac{0-L \bullet Y^p}{\|L\|^2}) L$
		\STATE $Z_{n+2m+3}^{p+1} \leftarrow U^p \max \{ \Gamma^{(p)},0\} U^{p T}$, where $U^P$ (resp. $\Gamma^p$) are the eigenvectors (resp. the eigenvalues) of $Y^p$
		\STATE $B'^{p+1} \leftarrow\frac{1}{n+2m+3} \Sigma_{i=1}^{n+2m+3} Z_i^{p+1}$
		 \IF {$\|B'^{p+1}-Z^p\|^2< \alpha $} 
		    \STATE stop $\leftarrow$ true
		 \ELSE
		    \STATE $B^{p+1} \leftarrow\beta^{p+1} B'^{p+1}+(1-\beta^{p+1})Z^p$ with $\beta^{p+1} \leftarrow\frac{\Sigma_{i=1}^{n+2m+3} \|Z_i^{p+1}-Z^p\|^2}{(n+2m+3)\|B'^{p+1}-Z^p\|^2}$
		    \STATE $Z^{p+1} \leftarrow\begin{cases}
    Z^p+\lambda(B^{p+1}-Z^p), & \text{if } p+1 \equiv k (\textrm{mod } k).\\
    B^{p+1} & \text{otherwise}.
  \end{cases}$
  \STATE $p++$
		 \ENDIF
  
  \ENDWHILE
  \STATE \textbf{return} $Z^p$

\end{algorithmic}
\label{algo2}
\end{algorithm}
{Solving} Problem~\eqref{bidual_y} using Algorithm~\ref{algo2} {requires the storage of} the matrices $M$, $O_j$ $j=1, \ldots, n$, $C_i$, $i=1, \ldots, m$, $Q$, $D_i$, $i=2, \ldots, m+1$ and $L$, that is in total $n+2m+4$ matrices of dimension $(2m+4) \times (2m+4)$. Nevertheless, there is a way to avoid storing these matrices, since Algorithm~\ref{algo2} {does not require} the whole matrices, but rather the result of operations that mostly include dot products of sparse matrices. {Then, doing} the calculations and {giving} the results needed in Lines 5, 7, 10, 12, 14 and 16 of Algorithm~\ref{algo2} in function of $A$, $b$, $\mu$, and $\Sigma$ {implies that} there will be no need for the matrices themselves. All these calculations are detailed in~\ref{sparse_computations}. This aspect is one of the contributions of this paper.

%

\section{Experimental results}\label{sec:experiments}
The experimental results aim at evaluating numerically the quality of the proposed solution by DFW Algorithm. As mentioned before, two ways of evaluating the quality of the solutions are possible: the first one is to compare with the exact solution proposed by CPLEX when solving optimally the BSOCP formulation of the problem. The other method is to compute an optimality gap obtained by the bidualization of the problem. For this, an important observation is that the bidual problem~\eqref{bidual_y} is an SDP problem. In order to compute this optimality gap, both CVXPY SDP solver and Pierra's algorithm {are used}.

First, the quality of the solution of DFW Algorithm {is evaluated} by the two methods mentioned before. Then, for the SDP relaxation, solutions obtained by CVXPY and Pierra's algorithm {are compared}, and the storage economy resulting from using Pierra's algorithm {is shown}. This storage economy is especially due to taking advantage of the matrices sparsity in Problem~\eqref{bidual_y}.

\subsection{Experimental setup}

The robust counterpart of the shortest path problem with an undirected grid graph {is considered} for different sizes. For a grid graph $L \times L$, the number of nodes is $n=L^2$, and the number of edges is $m=4 L (L -1)$. 
For the definition of Problem~\eqref{elli_prob}, the random mean vector $\mu$ and the random covariance matrix $\Sigma$ are chosen randomly, and {$\Omega$ is set to $1$}.

The implementation of both the computation of DFW robust solutions and the CVXPY based solver are written using {\tt{Python 3.8.5}} and Pierra's algorithm is implemented using {\tt{Matlab R2018b}}.

\subsection{Numerical evaluation of the heuristic approach DFW}

{This part contains comparisons between} the solutions obtained using DFW algorithm and CPLEX, {as well as computations of} the lower bound $d^{**}$ for the solution of DFW algorithm that is the solution of the SDP relaxation of the original problem. In this part, this lower bound is computed using CVXPY. Finally, the relative biduality gap $RBG$ (defined in~\eqref{rbg}) {is given,} which allows to evaluate the quality of the solution of DFW algorithm, and a performance ratio useful for comparison with other work.

For experiments with DFW algorithm, constant parameters are $\varepsilon=10^{-6}$ and $K=1000$. Table~\ref{lower_bound} shows results for problem sizes $L\in\{3, 4, \ldots, 10\}$. First of all, note that in all the cases processed, DFW algorithm gives the same solution as CPLEX. Second, concerning the relative gap, since we theoretically only have the weak duality for the {solved} problem, the biduality gap is not necessarily zero even if the solution is optimal. Thus, if the gap is small, it means that the heuristic solution is close to the optimal solution, but the opposite {may not be true}. Indeed a large gap does not mean that the heuristic solution is far from the optimal solution. In all the processed cases, the lower bound $d^{**}$ is less than the optimal solution $p^*$, which validates the developments and the computations. Then, the relative gap $RBG$ is between $0.1917$ and $0.3178$. This gap is a metric that allows to measure how far the heuristic approach is from the optimal solution. In other words, the heuristic solution is between $19.17\%$ and $31.78\%$ from optimality, in the worst case. To analyze more precisely the obtained gap and to compare with other work, a column is added in Table~\ref{lower_bound} with a metric named \textit{performance ratio} used in~\cite{karloff1999good} for the \textit{Max-Cut} problem. The \textit{performance ratio} $\rho$ has the following definition:

\begin{equation}
\rho = \frac{d^{**}}{g(\hat{x})}.
\end{equation}
$\rho$ is then the proportion of $d^{**}$ in $g(\hat{x})$. In {the considered} case, $0.7588\leq\rho\leq 0.8391$. This is comparable to $0.87$, the highest performance ratio obtained for the \textit{Max-Cut} problem.

It would be interesting to test other cases of larger problems where the comparison with CPLEX is not possible, and check if the relative gap stays in the same interval as the processed cases. Indeed, in the processed cases, DFW algorithm gives the optimal solution, and thus the gap only comes from $p^*-d^{**}$ (see Equation~\eqref{the_used_inequality}). 

Now that the evaluation of solutions given by DFW algorithm is done using both CPLEX and the relative gap computed by CVXPY, an issue remains as discussed in Section~\ref{subsection:SDPsolving}: CVXPY needs a huge amount of memory to store matrices. That has justified the use of an alternative approach with Pierra's algorithm using sparse computations detailed in~\ref{sparse_computations}. In the next section, numerical results obtained using Pierra's algorithm are presented, as well as the resulting gain in memory storage.

\begin{table}
\centering
\begin{tabular}{|p{1cm}||p{2.2cm}|p{2cm}|p{2cm}|p{2cm}|p{2.1cm}|}
 \hline
 $L$ & Solution of DFW $g(\hat{x})$ & Optimal solution by CPLEX $p^*$ & Lower bound by CVXPY $d^{**}$ & Relative gap RBG & Performance ratio $\rho$\\
 \hline
 3 & 223.8807 & 223.8807 & 169.8902 & 0.3178 & 0.7588\\
 \hline
 4 & 302.9097 & 302.9097 & 230.64099 & 0.3133 & 0.7614\\
 \hline
 5 & 381.3647 & 381.3647 & 292.6109 & 0.3033 & 0.7673\\
 \hline
 6 & 498.444952 & 498.444952 & 401.92866 & 0.2401 & 0.8064\\
 \hline
 7 & 524.41995 & 524.41995 & 422.3119 & 0.2418 & 0.8053\\
 \hline
 8 & 625.46595 & 625.46595 & 524.83906 & 0.1917 & 0.8391\\
 \hline
 9 & 659.0601 & 659.0601 & 542.6984 & 0.2144 & 0.8234\\
 \hline
 10 & 604.0187 & 604.0187 & 492.4042 & 0.2267 & 0.8152\\
 \hline
\end{tabular}
\caption{Comparison of the proposed solution by DFW with the optimal solution by CPLEX, and the lower bound by CVXPY.}
\label{lower_bound}
\end{table}

\subsection{Numerical results of Pierra's algorithm}

\begin{figure}
 \centering
 \includegraphics[scale=1]{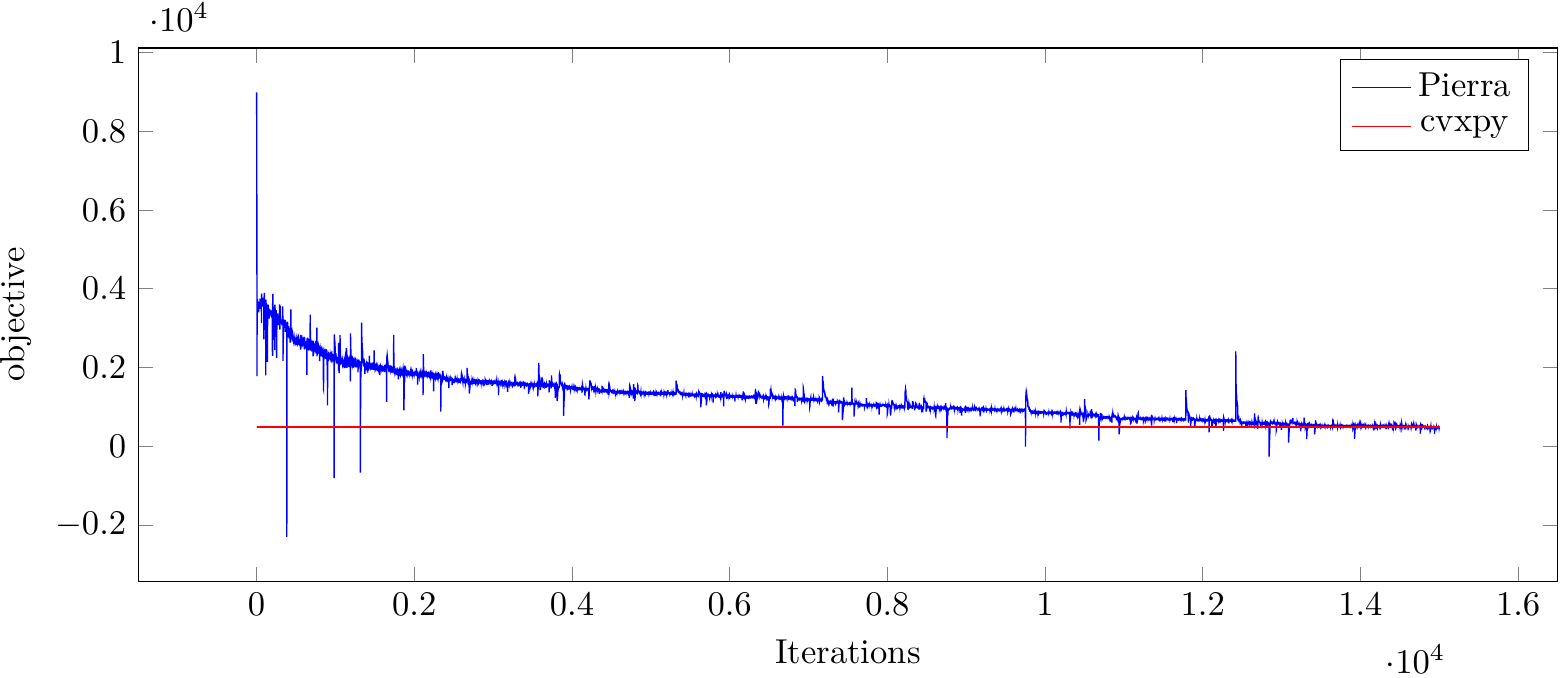}
 \caption{Evolution of the objective function along $15000$ iterations in Pierra's Algorithm compared to CVXPY's solution}
 \label{10_15000}
 \end{figure}

{This part shows} the results of Pierra's algorithm for Problem~\eqref{bidual_y} in comparison with the solution of CVXPY for problem sizes $L\in\{3, 4, \dots, 10\}$. For these experiments, constant parameters are $\varepsilon=1e-4$, $\lambda=0.5$, $k=3$ and $\alpha=10^{-8}$. Table~\ref{comparison_table} shows the computation time and memory storage needed for both the computation using CVXPY and of Pierra's algorithm, as well as the percentage of optimality of Pierra's solution compared to CVXPY after about 10000 iterations. The memory space saving is important. For $L=10$, {the proposed algorithm} {reduced the memory consumption} from 3.45 GigaBytes to 26 MegaBytes: a factor of $100$. In a reasonable computation time, that is however longer than the computation time of CVXPY, {Pierra's algorithm achieves} great percentages from optimality. Figure~\ref{10_15000} shows an example of the evolution of the objective function along the iterations of Pierra's algorithm for the problem size $L= 10$, compared to the optimal solution obtained by CVXPY. In this example $P=15000$, and $\varepsilon=10^{-4}$. A very good convergence can be observed at the last iterations shown in Table~\ref{comparison_table}: $99.93\%$ from optimality.

\begin{table}
\centering
\begin{tabular}{ |p{1cm}||p{1.3cm}|p{1.3cm}|p{1.7cm}|p{1.7cm}|p{4cm}|}
 \hline
  & \multicolumn{2}{c|}{Time(s)} & \multicolumn{2}{c|}{Storage needed (mB)} & \\
 \hline
 L & CVXPY & Pierra & CVXPY & Pierra & Optimality percentage of Pierra (\% CVXPY) \\
 \hline
 3 & 11 & 3.7 & 1.29792  & 0.13632 &  96.4\% \\
 4 & 49.6761 & 97.2 & 9.2   & 0.50496 & 77\% \\
 5 & 145.93 & 631 & 40.45   & 1.358848 & 86\%\\
 6 & 394.2456 & 1005.4 & 132.88    & 3.008448 & 92.2\%\\
 7 & 935.8 & 2275  & 358.82  & 5.841792 & 92.4\%\\
 8 & 2274.85 & 7826 & 841.73    & 10.32448 & 96\%\\
 9 & 4724.6 & 22338 & 1776.192  & 16.99968 & 97\%\\
 10 & 9244.87 & 63585 & 3451.17   & 26.488128 & 99.93\%\\
 \hline
\end{tabular}
\caption{Comparison between CVXPY and Pierra}
\label{comparison_table}
\end{table}

\subsection{Discussion}

In conclusion of the numerical experiments, it is possible to make the following comments. The lower bounds using Pierra's algorithm have been provided for small problem sizes. Thus, the contribution of this work is to propose a method to evaluate the solution of a heuristic algorithm for {Problem}~\eqref{elli_prob} without comparing it with CPLEX, but rather with a lower bound. For this, a challenge has been encountered, since it is well known that using the duality makes the problems easier but bigger, as the dual problem is usually polynomial, but has more variables and more constraints. This challenge has been tackled using Pierra's algorithm {with} {the} sparse computations. Here, the goal is twofold: first, put the algorithm proposed by Pierra in 1984 back in the spotlight for its efficiency even if it has not been used much. The second goal is to show the power of having an explicit algorithm instead of a black box solver. {Doing this has made the sparse computations possible}, reducing drastically the {memory} storage necessity.

Interesting future works involve going further in the problem sizes: starting from a grid of size $L=40$, the problem becomes computationally demanding, as CPLEX becomes unable of giving a solution, and CVXPY for the lower bound necessitates terabytes of {memory} storage. But before being able to realize {that}, some challenges concerning Pierra's algorithm should be dealt with, such as the stopping criteria on line 19 of Algorithm~\ref{algo2} and the performance of the algorithm {that has to be sped up}. One should note that the architecture of the algorithm allows a very easy parallelization, since the projections on each constraint space are independent (lines 6 to 17 of Algorithm~\ref{algo2}). Thus, {a parallel implementation} could speed up the algorithm.

%

\section{Conclusion}
\label{sec:conclusion}

This paper studies the Robust counterpart of the Shortest Path Problem (RSPP) in the case of correlated ellipsoidal uncertainty set. {This problem is NP-hard, and exact methods exist to solve it, such as BSOCP solvers. Moreover, a heuristic algorithm named DFW has been proposed in~\cite{al2020frank}}. More precisely, this work proposes a lower bound to validate heuristic approaches that solve the RSPP, such as DFW Algorithm. This lower bound computation replaces the comparison with exact solvers as a validation method. To compute the proposed lower bound, recall that it is the solution of an SDP problem that can be solved by CVXPY using interior-point methods. Unfortunately, the bidual problem is a big problem with much more constraints and more variables than the original problem. Thus, despite its polynomial nature, the resolution of this bidual problem is {very time} consuming and {needs a huge memory space}. Therefore, the sparsity of the matrices that define the problem has been exploited to replace the classical solver by a sparse version of Pierra's decomposition through formalization in a product space algorithm. All this is numerically tested, showing that, {due to} the results of this paper, a polynomial time evaluation of the quality of the solution of DFW heuristic is possible without having the memory storage issue of the bidual problem.

\section*{Acknowledgment}
\noindent This work has been supported by
the EIPHI Graduate school (contract "ANR-17-EURE-0002"). 
{Computations have been performed on the supercomputer facilities of Mésocenter de calcul de Franche-Comté in Besançon, France.}

\appendix
\section{Sparse computations}
\label{sparse_computations}
The aim of this appendix is to detail the computations needed in Algorithm~\ref{algo2}, and the replacements done to avoid the storage of the matrices $M$, $O_j$ $j=1, \ldots, n$, $C_i$, $i=1, \ldots, m$, $Q$, $D_i$, $i=2, \ldots, m+1$ and $L$. Recall that doing this enables us to express all the formulas in function only of $A$, $b$, $\mu$, and $\Sigma$, and thus to avoid the storage of $n+2m+4$ matrices {of dimension $2m+4 \times 2m+4$.}
\bigskip

\textbf{The operation in Line 5}
\begin{algorithmic}[1]
\STATE $Y^p=Z^p-\frac{\varepsilon}{2(n+2m+3)}M $
\end{algorithmic}
\textbf{can be replaced by }
\begin{algorithmic}[1]
\STATE $Y^p=Z^p$
\STATE $Y^p{}_{[1,2 \rightarrow m+1]} = Y^p{}_{[1,2 \rightarrow m+1]} - \frac{\varepsilon}{4(n+2m+3)} \mu^T $
\STATE $Y^p{}_{[2 \rightarrow m+1,1]} = Y^p{}_{[2 \rightarrow m+1,1]} - \frac{\varepsilon}{4(n+2m+3)} \mu $
\STATE $Y^p{}_{[1,2m+2]}=Y^p{}_{[1,2m+2]}-\frac{\varepsilon}{4(n+2m+3)}$
\STATE $Y^p{}_{[2m+2,1]}=Y^p{}_{[2m+2,1]}-\frac{\varepsilon}{4(n+2m+3)}$
\end{algorithmic}

\bigskip

\textbf{The operation in Line 7}
\begin{algorithmic}[1]
\STATE $Z_{j}^{p+1}=Y^p + (\frac{b_j-O_j \bullet Y^p}{\|O_j\|^2}) O_j$
\end{algorithmic}
\textbf{can be replaced by }
\begin{algorithmic}[1]
\STATE $Z_j^{p+1}=Y^p$
\STATE $Z_j^{p+1}[1,2:m+1] = Z_j^{p+1}[1,2:m+1]+\frac{a_j}{2} A_{j*}$
\STATE $Z_j^{p+1}[2:m+1,1] = Z_j^{p+1}[2:m+1,1]+\frac{a_j}{2} A_{j*}$
\end{algorithmic}
with $A_{j*}$ is the vector containing the $j$-th lign of $A$, $a_j=\frac{b_j-O_j \bullet Y^p}{\|O_j\|^2}=\frac{2b_j-\Sigma_{i=1}^m A_{ji} (Y^p_{i+1 \text{  }1}+Y^p_{1 \text{  }i+1})}{\Sigma_{i=1}^m A_{ji}^2}$, since
$\|O_j\|^2= \frac{1}{2} \Sigma_{i=1}^m A_{ji}^2$, and $O_j \bullet Y^p= \Sigma_{i=1}^m A_{ji} \frac{(Y^p_{i+1 \text{  }1}+Y^p_{1 \text{  }i+1})}{2}$.

\bigskip

\textbf{The operation in Line 10}
\begin{algorithmic}[1]
\STATE $Z_{n+1+i}^{p+1}=Y^p + (\frac{-C_i \bullet Y^p}{\|C_i\|^2}) C_i$
\end{algorithmic}
\textbf{can be replaced by}
\begin{algorithmic}[1]
\STATE $Z_{n+1+i}^{p+1}=Y^p$
\STATE $Z_{n+1+i}^{p+1}[1,2:m+1]=Z_{n+1+i}^{p+1}[1,2:m+1]-\frac{c_i}{2} (\Sigma^{\frac{1}{2} T})_i $
\STATE $Z_{n+1+i}^{p+1}[1,m+1+i]=Z_{n+1+i}^{p+1}[1,m+1+i]+\frac{c_i}{2}$
\STATE $Z_{n+1+i}^{p+1}[2:m+1,1]=Z_{n+1+i}^{p+1}[2:m+1,1]-\frac{c_i}{2} (\Sigma^{\frac{1}{2} T})_i $
\STATE $Z_{n+1+i}^{p+1}[m+1+i,1]=Z_{n+1+i}^{p+1}[m+1+i,1]+\frac{c_i}{2}$
\end{algorithmic}
with $c_i=\frac{-C_i \bullet Y^p}{\|C_i\|^2}=\frac{ \Sigma_{k=1}^m (\Sigma^{\frac{1}{2} T})_{ik} (Y^p_{k+1 \text{  }1}+Y^p_{1 \text{  }k+1})-Y^p_{m+i+1 \text{  }1}-Y^p_{1 \text{  }m+i+1}   }{1+\Sigma_{k=1}^{m} (\Sigma^{\frac{1}{2} T})_{ik}^2 }$, since $\|C_i\|^2=\frac{1}{2}(1+\Sigma_{k=1}^{m} (\Sigma^{\frac{1}{2} T})_{ik}^2)$, and $C_i \bullet Y^p=-\Sigma_{k=1}^m (\Sigma^{\frac{1}{2} T})_{ik} \frac{Y^p_{k+1 \text{  }1}+Y^p_{1 \text{  }k+1}}{2}+\frac{Y^p_{m+i+1 \text{  }1}+Y^p_{1 \text{  }m+i+1}}{2} $ 

\bigskip

\textbf{The operation in Line 12}
\begin{algorithmic}[1]
\STATE $Z_{n+m+2}^{p+1}=Y^p+(\frac{0-Q \bullet Y^p}{\|Q\|^2}) Q$
\end{algorithmic}
\textbf{can be replaced by}
\begin{algorithmic}[1]
\STATE $Z_{n+m+2}^{p+1}=Y^p$
\STATE $Z_{n+m+2}^{p+1}[m+1+i,m+1+i]=Z_{n+m+2}^{p+1}[m+1+i,m+1+i] + q $ for $i$ between $1$ and $m$.
\STATE $Z_{n+m+2}^{p+1}[2m+2,2m+2]=Z_{n+m+2}^{p+1}[2m+2,2m+2] - q $
\STATE $Z_{n+m+2}^{p+1}[2m+3,2m+3]=Z_{n+m+2}^{p+1}[2m+3,2m+3] + q $
\end{algorithmic}
with $q=\frac{0-Q \bullet Y^p}{\|Q\|^2}=\frac{\Sigma_{k=m+2}^{2m+1} Y^p_{kk}-Y^p_{2m+2 \text{ } 2m+2}+Y^p_{2m+3 \text{ } 2m+3} }{m+2}$, since $\|Q\|^2= m+2$, and $Q \bullet Y^p=-\Sigma_{k=m+2}^{2m+1} Y^p_{kk}+Y^p_{2m+2 \text{ } 2m+2}-Y^p_{2m+3 \text{ } 2m+3} $

\bigskip

\textbf{The operation in Line 14}
\begin{algorithmic}[1]
\STATE $Z_{n+m+1+i}^{p+1}=Y^p + (\frac{0-D_i \bullet Y^p}{\|D_i\|^2}) D_i$
\end{algorithmic}
\textbf{can be replaced by}
\begin{algorithmic}[1]
\STATE $Z_{n+m+1+i}^{p+1}=Y^p$
\STATE $Z_{n+m+1+i}^{p+1}[i,i]=Z_{n+m+1+i}^{p+1}[i,i] + d_i $
\STATE $Z_{n+m+1+i}^{p+1}[1,i]=Z_{n+m+1+i}^{p+1}[1,i] - \frac{d_i}{2} $
\STATE $Z_{n+m+1+i}^{p+1}[i,1]=Z_{n+m+1+i}^{p+1}[i,1] - \frac{d_i}{2} $
\end{algorithmic}
with $d_i=\frac{0-D_i \bullet Y^p}{\|D_i\|^2}=-\frac{2}{3}(Y^p[i,i]-\frac{Y^p[i,1]+Y^p[1,i]}{2}) $, since $\|D_i\|^2=\frac{3}{2}$, and $D_i \bullet Y^p=Y^p[i,i]-\frac{Y^p[i,1]+Y^p[1,i]}{2}$

\bigskip

\textbf{The operation in Line 16}
\begin{algorithmic}[1]
\STATE $Z_{n+2m+3}^{p+1}=Y^p+(\frac{0-L \bullet Y^p}{\|L\|^2}) L$
\end{algorithmic}
\textbf{can be replaced by}
\begin{algorithmic}[1]
\STATE $Z_{n+2m+3}^{p+1}=Y^p$
\STATE $Z_{n+2m+3}^{p+1}[2m+4,2m+4]=Z_{n+2m+3}^{p+1}[2m+4,2m+4] - l $
\STATE $Z_{n+2m+3}^{p+1}[1,2m+2]=Z_{n+2m+3}^{p+1}[1,2m+2] + \frac{l}{2} $
\STATE $Z_{n+2m+3}^{p+1}[2m+2,1]=Z_{n+2m+3}^{p+1}[2m+2,1] + \frac{l}{2} $
\end{algorithmic}
with $l=\frac{0-L \bullet Y^p}{\|L\|^2}=\frac{2}{3}( Y^p_{2m+4 \text{ } 2m+4}- \frac{Y^p_{2m+2 \text{ }1}+Y^p_{1 \text{ }2m+2} }{2}) $ since $\|L\|^2=\frac{3}{2}$ and $L \bullet Y^p= - Y^p_{2m+4 \text{ } 2m+4}+ \frac{Y^p_{2m+2 \text{ }1}+Y^p_{1 \text{ }2m+2} }{2}$

\bibliographystyle{elsarticle-num}
\bibliography{sample.bib}

\end{document}